\documentclass[12pt,thmsa]{article}
\usepackage{amsmath}
\usepackage{graphicx}
\usepackage{amsfonts}
\usepackage{amssymb}
\usepackage{mathrsfs}

\newtheorem{thm}{Theorem}

\newtheorem{rmk}{Remark}

\begin{document}

\begin{center}
{\Large \textbf{$k$-Sample problem based on generalized maximum mean discrepancy}}

\bigskip

Armando Sosthene Kali  BALOGOUN\textsuperscript{a} , Guy Martial  NKIET\textsuperscript{b}  and Carlos OGOUYANDJOU\textsuperscript{a}

\bigskip

\textsuperscript{a}Institut de Math\'ematiques et de Sciences Physiques, Porto Novo, B\'enin.
\textsuperscript{b}Universit\'{e} des Sciences et Techniques de Masuku,  Franceville, Gabon.

\bigskip

E-mail adresses : sosthene.balogoun@imsp-uac.org,  guymartial.nkiet@mathsinfo.univ-masuku.com,  ogouyandjou@imsp-uac.org.

\bigskip
\end{center}

\noindent\textbf{Abstract.}In this paper we deal with the problem of testing for the quality of $k$ probability distributions. We introduce a generalization of the maximum mean discrepancy that permits to characterize the null hypothesis. Then, an estimator of it is proposed as test statistic,  and its asymptotic distribution under the null hypothesis is derived. Simulations show that the introduced procedure outperforms classical ones.

\bigskip


\section{Introduction}
\label{Intro}
\noindent An important problem in statistics consists in testing whether two or more probability distributions are identical against the alternative that at least two of them   may differ. The case of two distributions, namely the two-sample problem, have been extensively studied, and there  are traditional approaches for dealing with it such as the Kolomogorv-Smirnov, Cram\'er-von Mises,  Anderson-Darling tests (see, e.g., Conover, 1980; Gibbons and Chakraborti, 2003)  and other nonparametric procedures like the Kruskal-Wallis test (Kruskal and Wallis, 1952). More recently, Gretton et al (2012) tackled this problem by using a kernel-based approach. They introduced the maximum mean discrepancy (MMD) and proposed an approach for the two-sample problem based on a test statistic which is  an unbiased estimator of this MMD. The interest of their approach is that it relies on kernel-based methods that allow the ability to work with high-dimensional and structured data (e.g., Harchaoui et al, 2013). The extension of procedures for the two-sample problem to the  case of more than two distributions, namely the $k$-sample problem, has been of great interest in the literature. In this vein, some of  the aforementioned traditional tests have been extended for dealing with the $k$-sample problem ($k\geq 2$). This is the case for the Komogorv-Smirnov test (Kiefer, 1959; Wolf and Naus, 1973), the  Cram\'er-von Mises test (Kiefer, 1959) and the Anderson-Darling test (Scholz and Stephens, 1987). More recently, Zhand and Wu (2007) introduced procedures based on the likelihood ratio and showed that their proposal leads to  more powerful tests that the traditionnal ones.

In this paper, we deal with the $k$-sample problem by extending the kernel-based approach of Gretton et al (2012). After recalling, in Section 2,  some facts about the reproducing kernel Hilbert spaces, we introduce in Section 3 a generalized maximum mean discrepancy (GMMD) that permits to characterize the null hypothesis related to the $k$-sample problem. Then, an estimator of this GMMD is introduced, in Section 4,  as test statistic, and its asymptotic distribution under the null  hypothesis is derived. Section 5, is devoted to the presentation of simulations made in order to evaluate performance of our proposal and to compare it with known  methods.

\section{Notation and preliminaries}\label{sec2}
\noindent In this section, we briefly recall the notion of reproducing kernel hilbert space (RKHS) and we define some elements related to it that are useful in this paper. Then, we specify the $k$-sample problem that we deal with.  

\bigskip

\noindent Letting  $(\mathcal{X},\mathcal{B})$ be a measurable space, where $\mathcal{X}$ is a topological space and $\mathcal{B}$ is the corresponding Borel $\sigma$-field, we consider a Hilbert space $\mathcal{H}$ of functions from $\mathcal{X}$ to $\mathbb{R}$, endowed with an inner product $<\cdot,\cdot>_\mathcal{H}$. This space is said to be a RKHS   if there exists a kernel, that is a symmetric positive semi-definite function $K:\mathcal{X}^2\rightarrow \mathbb{R} $, such that for any $f\in \mathcal{H}$ and any $x\in\mathcal{X}$, one has $K(x,\cdot )\in \mathcal{H}$ and  $f(x)=<f,K(x,\cdot)>_\mathcal{H}$. When $\mathcal{H}$ is a RKHS with kernel $K$,  the map $\Phi\,:\,x\in\mathcal{X}\mapsto K(x,\cdot)\in \mathcal{H}$ characterizes $K$ since one has
\[
K(x,y)=<\Phi(x),\Phi(y)>_{\mathcal{H}}
\]
for any $(x,y)\in\mathcal{X}^2$. It is called the feature map and it is an important tool when dealing with kernel methods for statistical problems.

\noindent Throughout this paper, we consider a RKHS $\mathcal{H}$ with kernel $K$ satisfying the following assumption:

\bigskip

$(\mathscr{A}_1):$  $\Vert  K\Vert_{\infty}:= \sup\limits_{(x,y)\in \mathcal{X}^2} K(x,y)<+ \infty .$\\

\bigskip

\noindent Now, let us consider a random variable $X$ taking values in $\mathcal{X}$ and with probability distribution $\mathbb{P}$. If $\mathbb{E}(\Vert \Phi(X)\Vert_\mathcal{H})=\int_{\mathcal{X}}\Vert\Phi(x)\Vert_\mathcal{H}d\mathbb{P}(x)<+\infty$, the mean element $m$ associated with  $X$ is defined for all functions $f\in \mathcal{H}$ as the unique element in  $\mathcal{H}$ satisfying,
\begin{eqnarray}
<m,f>_{\mathcal{H}}=\mathbb{E}\left(f(X)\right)=\int_{\mathcal{X}}f(x)d\mathbb{P}(x).
\end{eqnarray}
It is very important to note that if hypothesis  $(\mathscr{A}_1)$ is satisfied, then the mean element $m$  is well-defined. Its empirical counterpart, obtained from a i.i.d.  sample ${X_1,\cdots,X_n}$    of   $X$,  is given by:

\begin{equation}\label{empir}
\widehat{m}=\frac{1}{n} \sum_{i=1}^{n}K(X_i,\cdot).
\end{equation}
Now, for $k\in\mathbb{N}^\ast$ and $\ell=1,\cdots,k$, let us consider a random variable $X_\ell$ taking  values in $\mathcal{X}$ and with distribution $\mathbb{P}_\ell$. We are interested with the related  $k$-sample problem, that is the problem of testing for the hypothesis
\[
\mathscr{H}_0\,:\,\mathbb{P}_1=\mathbb{P}_2=\cdots=\mathbb{P}_k
\,\,\,\textrm{ 
against }\,\,\,
\mathscr{H}_1\,:\,\,\exists \,(j,\ell),\,\,\mathbb{P}_j\neq\mathbb{P}_\ell.
\]
For dealing with that problem, we will first introduce the notion of generalized maximum mean discrepancy which generalizes the maximum mean discrepancy of Gretton et al (2012).
\section{The generalized maximum mean discrepancy}
\noindent When dealing with the two-sample problem, Gretton et al (2012) introduced the maximum mean discrepancy (MMD), that is the largest difference in expectation over functions in the unit ball of a RKHS. More precisely, considering the unit ball  $\mathscr{F}=\{f\in\mathcal{H},\,\,\Vert f \Vert_\mathcal{H}\leq 1\}$ and $(\mathbb{P},\mathbb{Q})$ a pair of Borel probability mùeasures on $(\mathcal{X},\mathcal{B})$, they define the MMD as
\begin{equation*}
\textrm{MMD}(\mathscr{F},\mathbb{P},\mathbb{Q}):=\sup_{f\in\mathscr{F}}\bigg(\mathbb{E}_\mathbb{P}\left( f(X)\right)-\mathbb{E}_\mathbb{Q}\left( f(Y)\right)\bigg)
=\Vert m_1-m_2\Vert_\mathcal{H},
\end{equation*}
where $X$ (resp. $Y$) is a random variable with probability distribution $\mathbb{P}$ (resp. $\mathbb{Q}$) and $m_1$ (resp. $m_2$) is the mean element associated to  $X$ (resp. $Y$), and they show that $\textrm{MMD}(\mathscr{F},\mathbb{P},\mathbb{Q})=0$ if, and only if, $\mathbb{P}=\mathbb{Q}$. From this, we will introduce the  generalized maximum mean discrepancy (GMMD) which will characterize the hypothesis $\mathscr{H}_0$ of the $k$-sample problem described above.

\bigskip

\noindent\textbf{Definition 1.}
The generalized maximum mean discrepancy, related to  $\mathbb{P}_1, \cdots,\mathbb{P}_k$  and $\pi=\left(\pi_1,\cdots,\pi_k\right)\in (]0,1[)^k$  with $\sum_{\ell=1}^k\pi_\ell=1$, is the real denoted by $\textrm{GMMD}(\mathbb{P}_1, \cdots,\mathbb{P}_k;\pi)$ and  defined as
\begin{equation*}
 \textrm{GMMD}^2(\mathbb{P}_1, \cdots,\mathbb{P}_k;\pi)=\sum_{j=1}^{k}\sum_{\underset{\ell\neq j}{\ell=1}}^{k}\pi_\ell \,\textrm{MMD}^2(\mathscr{F},\mathbb{P}_j,\mathbb{P}_\ell)=\sum_{j=1}^{k}\sum_{\underset{\ell\neq j}{\ell=1}}^{k}\pi_\ell \parallel m_j-m_k\parallel_{\mathcal{H}}^2
\end{equation*}
where  $m_\ell$ is the mean element associated to $\mathbb{P}_\ell$.

\bigskip

\noindent The definition of MMD appears to be a particular case of the above definition  obtained for  $k=2$.  The hypothesis $\mathscr{H}_0$ can be characterized by means of the GMMD. Indeed, it is easy to check that   this hypothesis  is true if, and only if,   $\textrm{GMMD}(\mathbb{P}_1, \cdots,\mathbb{P}_k;\pi)=0$. Then, a test statistic for the $k$-sample problem may be chosen by taking an estimator of  $ \textrm{GMMD}^2(\mathbb{P}_1, \cdots,\mathbb{P}_k;\pi)$. 
\section{Test statistic and   asymptotic distribution}
\noindent For $j=1,2,\cdots,k$,  let $\{X_1^{(j)},\cdots,X_{n_j}^{(j)}\}$  be an i.i.d. sample in  $\mathcal{X}$ with commmon distribution  $\mathbb{P}_j$   with  mean element
 $m_j$.  Putting $n=\sum_{j=1}^kn_j$, we assume that the following condition holds:
\bigskip

\noindent $(\mathscr{A}_2):$  For $j\in\{1,\cdots,k\}$, one has $\lim_{n_j\rightarrow +\infty}\frac{n_j}{n}=\rho_j$, where $\rho_j$ is a real belonging to $]0,1[$.

\bigskip

\noindent  We will first introduce a test statistic for the $k$-sample problem by taking an estimator of  $\textrm{GMMD}^2(\mathbb{P}_1, \cdots,\mathbb{P}_k;\pi)$ where  $\pi=\left(\rho_1,\cdots,\rho_k\right)$ , then we will derive its asymptotic distribution under $\mathscr{H}_0$.
\subsection{Test statistic}
\noindent We know from Lemma 6 in  Gretton et al (2012) that an unbiased estimator of  $\textrm{MMD}^2(\mathscr{F},\mathbb{P}_j,\mathbb{P}_\ell)$ is given by
\begin{eqnarray}\label{statdev}
 \widehat{\Gamma}_{j\ell}^{(n)}&=&\sum_{i=1}^{n_j} \sum_{\underset{r\neq i}{r=1}}^{n_j}  \frac{1}{n_j(n_j-1)}K(X_i^{(j)},X_r^{(j)})+\sum_{i=1}^{n_\ell} \sum_{\underset{r\neq i}{r=1}}^{n_\ell} \frac{1}{n_\ell(n_\ell-1)}K(X_i^{(\ell)},X_r^{(\ell)})\nonumber\\
& &-2\sum_{i=1}^{n_j} \sum_{r=1}^{n_\ell}\frac{1}{n_jn_\ell}K(X_i^{(j)},X_r^{(\ell)}) .
\end{eqnarray}
Then,    we take as test statistic the  estimator of  $\textrm{GMMD}^2(\mathbb{P}_1, \cdots,\mathbb{P}_k;\pi)$ defined as
\begin{eqnarray*}
\widehat{T}_n=\sum_{j=1}^k \sum_{\underset{\ell\neq j}{\ell=1}}^k P_\ell\,\widehat{\Gamma}_{j\ell}^{(n)},
\end{eqnarray*}
where    $P_\ell=\frac{n_\ell}{n}$. 

\subsection{Asymptotic distribution of  $\widehat{T} _n$ under  $\mathscr{H}_0$}
\noindent Under $\mathscr{H}_0$, one has $m_1=\cdots=m_2=m$;   let us then  consider the centered kernel $\widetilde{K}$ defined as
	$
	\widetilde{K}(x,y)=<\Phi(x)-m,\Phi(y)-m>_{\mathcal{H}}
=K(x,y)-\mathbb{E}\left(K(X,x)\right)-\mathbb{E}\left(K(X,y)\right)+\mathbb{E}\left(K(X,X^\prime)\right)$,
where $X^\prime$ is a random variable independent of $X$ and with the same distribution $\mathbb{P}$, and the sequence $\{\lambda_p\}_{p\geq 1}$ of eigenvalues of the integral operator associated to $\widetilde{K}$. Then, we have :

\bigskip

\begin{thm}\label{loi}
We suppose that the assumptions  $(\mathscr{A}_1)$  and $(\mathscr{A}_2)$ hold. Then, under  $\mathscr{H}_0$,  as $\min_{1\leq j\leq k}(n_j)\rightarrow +\infty$, one has 
\begin{eqnarray}
n \widehat{T} _n\stackrel{\mathscr{D}}{\rightarrow}\sum_{p=1}^{+\infty}\lambda_p\left\{(k-2)(Z_p-k)+\sum_{j=1}^{k}\left\{\rho_j^{-1}(Y^2_{p,j}-1)-2\sum_{\underset{\ell\neq j}{\ell=1}}^k\rho_\ell^{1/2}\rho_j^{-1/2}Y_{p,j}Y_{p,\ell}\right\}\right\}
\end{eqnarray}
where $(Y_{p,j})_{p\geq 1,\,1\leq j\leq k}$ is a sequence of independent random variables having the $\mathcal{N}(0,1)$  distribution, and  $(Z_p)_{p\geq 1}$ is a sequence of independent random variables having the $\chi^2_k$ distribution.
\end{thm}
\noindent\textit{Proof.}
  Letting  $Z$ and $Z^{'}$ be two independent  random variables with distribution $\mathbb{P}=\mathbb{P}_\ell$, $\ell=1,\cdots,k$,  we have from (\ref{statdev}):
\begin{eqnarray*}\label{ac}
\widehat{\Gamma}_{j\ell}^{(n)}&=&\frac{1}{n_j(n_j-1)}\sum_{i=1}^{n_j} \sum_{\underset{r\neq i}{r=1}}^{n_j}\bigg\{\widetilde{K}(X_i^{(j)},X_r^{(j)})+\mathbb{E}\left(K(X^{(j)}_1,Z)\right)\\
& &\hspace{3.5cm}+\mathbb{E}\left(K(Z,X^{(\ell)}_1)\right)-\mathbb{E}\left(K(Z,Z')\right)\bigg\}\nonumber\\
& &+\frac{1}{n_\ell(n_\ell-1)}\sum_{i=1}^{n_\ell} \sum_{\underset{r\neq i}{r=1}}^{n_\ell}\bigg\{\widetilde{K}(X_i^{(\ell)},X_r^{(\ell)})+\mathbb{E}\left(K(X^{(j)}_1,Z)\right)\\
& &\hspace{3.5cm}+\mathbb{E}\left(K(Z,X^{(\ell)}_1)\right)-\mathbb{E}\left(K(Z,Z')\right)\bigg\}\nonumber\\
& &-\frac{2}{n_jn_\ell}\sum_{i=1}^{n_j}\sum_{r=1}^{n_\ell}\bigg\{\widetilde{K}(X_i^{(j)},X_r^{(\ell)})+\mathbb{E}\left(K(X^{(j)}_1,Z)\right)\\
& &\hspace{3.5cm}+\mathbb{E}\left(K(Z,X^{(\ell)}_1)\right)-\mathbb{E}\left(K(Z,Z')\right)\bigg\}\nonumber\\
&=&\frac{1}{n_j(n_j-1)}\sum_{i=1}^{n_j} \sum_{\underset{r\neq i}{r=1}}^{n_j}\widetilde{K}(X_i^{(j)},X_r^{(j)})\\
& &+\frac{1}{n_\ell(n_\ell-1)}\sum_{i=1}^{n_\ell} \sum_{\underset{r\neq i}{r=1}}^{n_\ell}\widetilde{K}(X_i^{(\ell)},X_r^{(\ell)})\nonumber\\
& &-\frac{2}{n_jn_\ell}\sum_{i=1}^{n_j}\sum_{r=1}^{n_\ell}\widetilde{K}(X_i^{(j)},X_r^{(\ell)})
\end{eqnarray*}
and,  from (\ref{statdev}),  $n \widehat{T} _n=A_n+B_n$, where
\begin{equation*}\label{an}
A_n=\sum_{j=1}^k \sum_{\underset{\ell\neq j}{\ell=1}}^k\bigg\{\sum_{i=1}^{n_j} \sum_{\underset{r\neq i}{r=1}}^{n_j} \frac{nP_\ell}{n_j(n_j-1)}\widetilde{K}(X_i^{(j)},X_r^{(j)})
+\sum_{i=1}^{n_\ell} \sum_{\underset{r\neq i}{r=1}}^{n_\ell}\frac{nP_\ell}{n_\ell(n_\ell-1)}\widetilde{K}(X_i^{(\ell)},X_r^{(\ell)})\bigg\}
\end{equation*}
and 
\begin{equation*}\label{cn}
B_n=-2\sum_{j=1}^k \sum_{\underset{\ell\neq j}{\ell=1}}^k\sum_{i=1}^{n_j} \sum_{r=1}^{n_\ell}\frac{nP_\ell}{n_jn_\ell}\widetilde{K}(X_i^{(j)},X_r^{(\ell)}).
\end{equation*}
Since $K$ is bounded (from assumption ($\mathscr{A}_1$)), the integral operator $S_{\widetilde{K}}$ associated to $\widetilde{K}$ is a Hibert-Schimdt operator and has, therefore, a system  $\{e_p\}_{p\geq 1}$ of eigenfunctions that is an orthonormal basis of $L^2(\mathbb{P})$. Thus,  $\widetilde{K}(x,y)=\sum_{p=1}^{+\infty}\lambda_p\,e_p(x)\,e_p(y)$ and
\begin{eqnarray}\label{an}
A_n&=&\sum_{j=1}^k \sum_{\underset{\ell\neq j}{\ell=1}}^k\bigg\{\sum_{i=1}^{n_j} \sum_{\underset{r\neq i}{r=1}}^{n_j}\sum_{p=1}^{+\infty}\frac{nP_\ell}{n_j(n_j-1)} \lambda_{p}e_{p}(X_i^{(j)})e_{p}(X_r^{(j)})\nonumber\\
& &+\sum_{i=1}^{n_\ell} \sum_{\underset{r\neq i}{r=1}}^{n_\ell}\sum_{p=1}^{+\infty}\frac{nP_\ell}{n_\ell(n_\ell-1)}\lambda_{p}e_{p}(X_i^{(\ell)})e_{p}(X_r^{(\ell)})\bigg\}\nonumber\\
&=&\sum_{j=1}^k \sum_{i=1}^{n_j} \sum_{\underset{r\neq i}{r=1}}^{n_j}\sum_{p=1}^{+\infty}\frac{n}{n_j(n_j-1)}\bigg(\sum_{\underset{\ell\neq j}{\ell=1}}^kP_\ell\bigg) \lambda_{p}e_{p}(X_i^{(j)})e_{p}(X_r^{(j)})\nonumber\\
& &+\sum_{\ell=1}^k\sum_{i=1}^{n_\ell} \sum_{\underset{r\neq i}{r=1}}^{n_\ell}\sum_{p=1}^{+\infty}\frac{n(k-1)P_\ell}{n_\ell(n_\ell-1)}\lambda_{p}e_{p}(X_i^{(\ell)})e_{p}(X_r^{(\ell)})\nonumber\\
&=&\sum_{j=1}^k \sum_{i=1}^{n_j} \sum_{\underset{r\neq i}{r=1}}^{n_j}\sum_{p=1}^{+\infty}\frac{(1-P_j)n}{n_j(n_j-1)} \lambda_{p}e_{p}(X_i^{(j)})e_{p}(X_r^{(j)})\nonumber\\
& &+\sum_{\ell=1}^k\sum_{i=1}^{n_\ell} \sum_{\underset{r\neq i}{r=1}}^{n_\ell}\sum_{p=1}^{+\infty}\frac{n(k-1)P_\ell}{n_\ell(n_\ell-1)}\lambda_{p}e_{p}(X_i^{(\ell)})e_{p}(X_r^{(\ell)})\nonumber\\
&=&\sum_{j=1}^{k}\left\{\sum_{p=1}^{+\infty}\sum_{i=1}^{n_j}\sum_{\underset{r\neq i}{r=1}}^{n_j}\frac{n[1+(k-2)P_j]}{n_j(n_j-1)} \lambda_{p}e_{p}(X_i^{(j)})e_{p}(X_r^{(j)})\right\}\nonumber\\
&=&\sum_{j=1}^{k}\frac{n[1+(k-2)P_j]}{n_j-1}\bigg\{\sum\limits_{p=1}^{+\infty}\lambda_{p}\bigg[\left( \frac{1}{\sqrt{n_{j}}}\sum\limits_{i=1}^{n_{j}}e_{p}(X^{(j)}_{i})\right)^2\nonumber\\
& &\hspace{5cm}-\frac{1}{n_{j}}\sum_{i=1}^{n_{j}}e^2_{p}(X^{(j)}_{i}) \bigg]\bigg\},
\end{eqnarray} 
\begin{eqnarray}\label{bn}
B_n&=&-2\sum_{j=1}^k \sum_{\underset{\ell\neq j}{\ell=1}}^k\sum_{i=1}^{n_j} \sum_{r=1}^{n_\ell}\sum_{p=1}^{+\infty}\frac{-2nP_\ell}{n_jn_\ell}\lambda_{p}e_{p}(X_i^{(\ell)})e_{p}(X_r^{(j)})\nonumber\\
&=&-2\sum_{p=1}^{\infty}\lambda_p\sum_{j=1}^{k}\left\{\frac{P_j^{-1/2}}{\sqrt{n_j}}\sum_{i=1}^{n_j}e_{p}(X_i^{(j)})\left(\sum_{\ell=1}^{k}\frac{P_\ell^{1/2}}{\sqrt{n_\ell}}\sum_{r=1}^{n_\ell}e_{p}(X_r^{(\ell)})-\frac{P_j^{1/2}}{\sqrt{n_j}}\sum_{r=1}^{n_j}e_p(X_r^{(j)})\right)\right\}\nonumber\\
&=&2\sum_{p=1}^{+\infty}\lambda_p\sum_{j=1}^{k}\left[\frac{1}{\sqrt{n_j}}\sum_{i=1}^{n_j}e_{p}(X_i^{(j)})\right]^2\nonumber\\
& &-2\sum_{p=1}^{+\infty}\lambda_p\left\{\left[\sum_{j=1}^{k}\frac{P_j^{-1/2}}{\sqrt{n_j}}\sum_{i=1}^{n_j}e_{p}(X_i^{(j)})\right]\left[\sum_{\ell=1}^{k}\frac{P_\ell^{1/2}}{\sqrt{n_\ell}}\sum_{r=1}^{n_\ell}e_{p}(X_r^{(\ell)})\right]\right\}\nonumber\\
&=&2\sum_{p=1}^{\infty}\lambda_p\sum_{j=1}^{k}\left[\frac{1}{\sqrt{n_j}}\sum_{i=1}^{n_j}e_{p}(X_i^{(j)})\right]^2+\sum_{p=1}^{\infty}\lambda_p\left[\sum_{j=1}^{k}\frac{1-P_j}{\sqrt{P_jn_j}}\sum_{i=1}^{n_j}e_{p}(X_i^{(j)})\right]^2\nonumber\\
&-&\sum_{p=1}^{\infty}\lambda_p\left\{\left[\sum_{j=1}^{k}\frac{P_j^{-1/2}}{\sqrt{n_j}}\sum_{i=1}^{n_j}e_{p}(X_i^{(j)})\right]^2+\left[\sum_{j=1}^{k}\frac{P_j^{1/2}}{\sqrt{n_j}}\sum_{r=1}^{n_j}e_{p}(X_r^{(j)})\right]^2\right\}.
 \end{eqnarray}
Then, from (\ref{an}) and (\ref{bn}),   $n \widehat{T} _n =\sum\limits_{p=1}^{+\infty}\lambda_p\,\mathcal{W}_{n,p}$, where $\mathcal{W}_{n,p}=\phi_n(\mathcal{U}_{n,p})-\psi_n(\mathcal{V}_{n,p})$ and
\[
\mathcal{U}_{n,p}=(U_{n_1,p},\cdots, U_{n_k,p})\,\,\, \,\,\,\,\,\,(\textrm{resp. }  \mathcal{V}_{n,p}=(V_{n_1,p},\cdots, V_{n_k,p}))
\] 
with $U_{n_j,p}=n_j^{-1/2}\sum_{i=1}^{n_j}e_{p}(X_i^{(j)})$ (resp. $V_{n_j,p}=n_j^{-1}\sum_{i=1}^{n_j}e_{p}^{2}(X_i^{(j)})$), the maps  $\phi_n$ and $\psi_n$   from $\mathbb{R}^p$ to $\mathbb{R}$ being  defined as
 \begin{eqnarray*}
 \phi_n(x)&=&\sum_{j=1}^{k}\bigg(\frac{n(1+(k-2)P_j)}{n_j-1}+
 2\bigg)x_j^2+\left(\sum_{j=1}^{k}\frac{1-P_j}{\sqrt{P_j}}x_j\right)^2\\
& &-\left(\sum_{j=1}^{k}P_j^{-1/2}x_j\right)^2-\left(\sum_{j=1}^{k}P_j^{1/2}x_j\right)^2
 \end{eqnarray*}
and
\[
\psi _n(x)=\sum_{j=1}^{k}\frac{n(1+(k-2)P_j)}{n_j-1}x_j.
\]
Since, for  $(j,\ell)\in \{1,2,\cdots,k\}^2$ with $j\neq \ell$, $U_{n_j,p}$ and   $U_{n_\ell,p}$ are independent and, from  the   central  limit theorem, 
$U_{n_j,p}\stackrel{\mathscr{D}}{\rightarrow}Y_{j,p}$ as $n_j\rightarrow +\infty$ where  $Y_{j,p}\sim \mathcal{N}(0,1)$, we deduce that   $\mathcal{U}_{n,p}\stackrel{\mathscr{D}}{\rightarrow}\mathcal{U}_{p}:=(Y_{1,p},\cdots,Y_{k,p})$ as $\min_{1\leq j\leq k}(n_j)\rightarrow +\infty$ where, for $j\neq \ell$,   $Y_{j,p}$ and $Y_{\ell,p}$ are independent variables having $\mathcal{N}(0,1)$ distribution. On the other hand, from law of large numbers, $\mathcal{V}_{n,p}$ converges in probability to ${\large \textrm{\textbf{1}}}_k:=(1,1,\cdots,1)$. Then, considering the maps  $\phi$ and $\psi$   from $\mathbb{R}^p$ to $\mathbb{R}$   defined as
 \begin{eqnarray*}
 \phi(x)&=&\sum_{j=1}^{k}\bigg(\rho_j^{-1}+(k-2)+
 2\bigg)x_j^2+\left(\sum_{j=1}^{k}\rho_j^{-1/2}\left(1-\rho_j\right)x_j\right)^2\\
& &-\left(\sum_{j=1}^{k}\rho_j^{-1/2}x_j\right)^2-\left(\sum_{j=1}^{k}\rho_j^{1/2}x_j\right)^2
 \end{eqnarray*}
and
\[
\psi(x)=\sum_{j=1}^{k}\left(\rho_j^{-1}+k-2\right)x_j,
\]
and putting  $\mathcal{W}_{p}=\phi(\mathcal{U}_p)-\psi({\large \textrm{\textbf{1}}}_k)$, we will show that  $n \widehat{T} _n \stackrel{\mathscr{D}}{\rightarrow}\sum\limits_{p=1}^{+\infty}\lambda_p\,\mathcal{W}_{p}$ as $\min_{1\leq j\leq k}(n_j)\rightarrow +\infty$. First, denoting by $\varphi_U$ the characteristic function of the random variable $U$,  putting $\widehat{S}_n=n \widehat{T} _n $ and $\widehat{S}_n^{(q)}=\sum\limits_{p=1}^{q}\lambda_p\,\mathcal{W}_{n,p}$ for $q\in\mathbb{N}^\ast$, and using the inequality $\vert e^{iz}-1\vert\leq\vert z\vert$ which holds for any $z\in\mathbb{R}$, we have for any $t\in\mathbb{R}$:
\begin{eqnarray*}
\left\vert\varphi_{\widehat{S}_n}(t)-\varphi_{\widehat{S}_n^{(q)}}(t)\right\vert&\leq& 
\mathbb{E}\left(\left\vert e^{it\widehat{S}_n}-e^{it\widehat{S}_n^{(q)}}\right\vert\right)
=\mathbb{E}\left(\left\vert e^{it\left(\widehat{S}_n-\widehat{S}_n^{(q)}\right)}-1\right\vert\right)\\
&\leq&\left\vert t\right\vert\mathbb{E}\left(\left\vert  \widehat{S}_n-\widehat{S}_n^{(q)} \right\vert\right)
\leq\left\vert t\right\vert\sum\limits_{p=q+1}^{+\infty}\lambda_p\mathbb{E}\left(\left\vert  \mathcal{W}_{n,p} \right\vert\right)
\end{eqnarray*}
and
\begin{eqnarray*}
\mathbb{E}\left(\left\vert  \mathcal{W}_{n,p} \right\vert\right)&=&\mathbb{E}\left(\left\vert  \phi_n(\mathcal{U}_{n,p})-\psi_n(\mathcal{V}_{n,p}) \right\vert\right)\\
&\leq&\sum_{j=1}^{k}\frac{n[1+(k-2)P_j]}{n_j-1}\mathbb{E}\left(U_{n_j,p}^2\right)+\mathbb{E}\left(\left[\sum_{j=1}^{k}\frac{1-P_j}{\sqrt{P_j}}U_{n_j,p}\right]^2\right)\\
& &+\mathbb{E}\left(\left[\sum_{j=1}^{k}P_j^{-1/2}U_{n_j,p}\right]^2\right)+\mathbb{E}\left(\left[\sum_{j=1}^{k}P_j^{1/2}U_{n_j,p}\right]^2\right)\\
& &+\sum_{j=1}^{k}\frac{n[1+(k-2)P_j]}{n_j-1}\mathbb{E}\left(V_{n_j,p}\right).
\end{eqnarray*}
Since $\mathbb{E}\left(e_p^2(X_i^{(j)})\right)=1$  and  $\mathbb{E}\left(e_p(X_i^{(j)})\,e_p(X_r^{(\ell)})\right)=\delta_{ir}\delta_{j\ell}$, it follows
\[
\mathbb{E}\left(V_{n_j,p}\right)=\frac{1}{n_j}\sum_{i=1}^{n_j}\mathbb{E}\left(e_p^2(X_i^{(j)})\right)=1
\]
and
\begin{eqnarray*}
\mathbb{E}\left(U_{n_j,p}^2\right)=\frac{1}{n_j}\sum_{i=1}^{n_j}\mathbb{E}\left(e_p^2(X_i^{(j)})\right)+\frac{1}{n_j}\sum_{i=1}^{n_j}\sum_{\underset{r\neq i}{r=1}}^{n_j}
\mathbb{E}\left(e_p(X_i^{(j)})\,e_p(X_r^{(j)})\right)=1,
\end{eqnarray*} 

\begin{eqnarray*}
\mathbb{E}\left(\left[\sum_{j=1}^{k}\frac{1-P_j}{\sqrt{P_j}}U_{n_j,p}\right]^2\right)&=&\sum_{j=1}^{k}\frac{(1-P_j)^2}{P_j}\mathbb{E}\left(U_{n_j,p}^2\right)\\
& &+\sum_{j=1}^{k}\sum_{\underset{\ell\neq j}{\ell=1}}^{k}\frac{(1-P_j)(1-P_\ell)}{\sqrt{P_jP_\ell}}\mathbb{E}\left(U_{n_j,p}U_{n_\ell,p}\right)\\
&=&\sum_{j=1}^{k}\frac{(1-P_j)^2}{P_j}\\
& &+\sum_{j=1}^{k}\sum_{\underset{\ell\neq j}{\ell=1}}^{k}\sum_{i=1}^{n_j}\sum_{r=1}^{n_\ell}\frac{(1-P_j)(1-P_\ell)}{\sqrt{n_jn_\ell P_jP_\ell}}\mathbb{E}\left(e_p(X_i^{(j)})\,e_p(X_r^{(\ell)})\right)\\
&=&\sum_{j=1}^{k}\frac{(1-P_j)^2}{P_j},
\end{eqnarray*} 
\begin{eqnarray*}
\mathbb{E}\left(\left[\sum_{j=1}^{k}P_j^{-1/2}U_{n_j,p}\right]^2\right)&=&\sum_{j=1}^{k} P_j^{-1} \mathbb{E}\left(U_{n_j,p}^2\right) +\sum_{j=1}^{k}\sum_{\underset{\ell\neq j}{\ell=1}}^{k} P_j^{-1} P_\ell^{-1}\mathbb{E}\left(U_{n_j,p}U_{n_\ell,p}\right)\\
&=&\sum_{j=1}^{k} P_j^{-1} +\sum_{j=1}^{k}\sum_{\underset{\ell\neq j}{\ell=1}}^{k}\sum_{i=1}^{n_j}\sum_{r=1}^{n_\ell} P_j^{-1} P_\ell^{-1}\mathbb{E}\left(e_p(X_i^{(j)}\,e_p(X_r^{(\ell)})\right)\\
&=&\sum_{j=1}^{k} P_j^{-1},
\end{eqnarray*}
\begin{eqnarray*}
\mathbb{E}\left(\left[\sum_{j=1}^{k}P_j^{1/2}U_{n_j,p}\right]^2\right)&=&\sum_{j=1}^{k} P_j\mathbb{E}\left(U_{n_j,p}^2\right) +\sum_{j=1}^{k}\sum_{\underset{\ell\neq j}{\ell=1}}^{k} P_j^{1/2} P_\ell^{1/2}\mathbb{E}\left(U_{n_j,p}U_{n_\ell,p}\right)\\
&=&\sum_{j=1}^{k} P_j +\sum_{j=1}^{k}\sum_{\underset{\ell\neq j}{\ell=1}}^{k}\sum_{i=1}^{n_j}\sum_{r=1}^{n_\ell} P_j^{1/2} P_\ell^{1/2}\mathbb{E}\left(e_p(X_i^{(j)}\,e_p(X_r^{(\ell)})\right)\\
&=&\sum_{j=1}^{k} P_j=1.
\end{eqnarray*}
Thus
\begin{eqnarray*}
\mathbb{E}\left(\left\vert  \mathcal{W}_{n,p} \right\vert\right)
&\leq&\sum_{j=1}^{k}\left(\frac{2n[1+(k-2)P_j]}{n_j-1}+\frac{(1-P_j)^2+1}{P_j}+1\right)
\end{eqnarray*}
and since 
\[
\lim_{n_j\rightarrow +\infty}\left(\frac{2n[1+(k-2)P_j]}{n_j-1}+\frac{(1-P_j)^2+1}{P_j}+1\right)=2(\rho_j^{-1}+k-2)+\frac{(1-\rho_j)^2+1}{\rho_j}+1,
\]
there exists $n_j^0\in\mathbb{N}^\ast$ such that, for $n_j\geq n_j^0$ ($j=1,\cdots, k$), we have
\begin{eqnarray*}
\mathbb{E}\left(\left\vert  \mathcal{W}_{n,p} \right\vert\right)
&\leq&\sum_{j=1}^{k}\left(2(\rho_j^{-1}+k-2)+\frac{(1-\rho_j)^2+1}{\rho_j}+2\right)
\end{eqnarray*}
and, therefore, 
\begin{eqnarray*}
\left\vert\varphi_{\widehat{S}_n}(t)-\varphi_{\widehat{S}_n^{(q)}}(t)\right\vert
\leq\left\vert t\right\vert\sum_{j=1}^{k}\left(2(\rho_j^{-1}+k-2)+\frac{(1-\rho_j)^2+1}{\rho_j}+2\right)\sum\limits_{p=q+1}^{+\infty}\lambda_p.
\end{eqnarray*}
Since $\sum_{p=1}^{+\infty}\lambda_p<+\infty$, we deduce  that $\lim_{q\rightarrow +\infty}\left\vert\varphi_{\widehat{S}_n}(t)-\varphi_{\widehat{S}_n^{(q)}}(t)\right\vert=0$. Then, for all $\varepsilon>0$  there exists $q_0\in\mathbb{N}^\ast$ such that 
\begin{equation}\label{eps1}
\left\vert\varphi_{\widehat{S}_n}(t)-\varphi_{\widehat{S}_n^{(q)}}(t)\right\vert<\frac{\varepsilon}{3}
\end{equation}
for $q\geq q_0$. Secondly, let us consider  $S_q=\sum_{p=1}^q\lambda_p\mathcal{W}_{p}$ and show that we have   $\widehat{S}_n^{(q)}\stackrel{\mathscr{D}}{\rightarrow}S_q$  as $\min_{1\leq j\leq k}(n_j)\rightarrow +\infty$. It suffices to prove that  $\widehat{S}_n^{(q)}-S_q$ converges in probability to $0$. For doing that we first  consider the  inequality
\begin{eqnarray}\label{ineq1}
\left\vert \widehat{S}_n^{(q)}-S_q\right\vert
&\leq&\sum_{p=1}^q\lambda_p\left\vert\mathcal{W}_{n,p}-\mathcal{W}_{p}    \right\vert\leq\sum_{p=1}^q\lambda_p\bigg(\left\vert\phi_n(\mathcal{U}_{n,p})-\phi(\mathcal{U}_p)\right\vert\nonumber\\
& &\hspace{4.5cm}+\left\vert \psi_n(\mathcal{V}_{n,p})-\psi({\large \textrm{\textbf{1}}}_k)\right\vert\bigg)\nonumber\\
&\leq &\sum_{p=1}^q\lambda_p\bigg(\left\vert\phi_n(\mathcal{U}_{n,p})-\phi(\mathcal{U}_{n,p})\right\vert+\left\vert\phi(\mathcal{U}_{n,p})-\phi(\mathcal{U}_p)\right\vert\nonumber\\
& &\hspace{1.5cm}+\left\vert \psi_n(\mathcal{V}_{n,p})- \psi(\mathcal{V}_{n,p})\right\vert+\left\vert  \psi(\mathcal{V}_{n,p})-\psi({\large \textrm{\textbf{1}}}_k)\right\vert\bigg).
\end{eqnarray}
Further,
\begin{eqnarray}\label{ineq2}
\left\vert\psi_n(\mathcal{V}_{n,p})-\psi(\mathcal{V}_{n,p})\right\vert\
&\leq& \sum_{j=1}^{k}\left\vert\frac{n(1+(k-2)P_j)}{n_j-1}-\rho^{-1}_j-k+2\right\vert\,\left\vert \mathcal{V}_{n,p}\right\vert
\end{eqnarray}
and, using $a^2-b^2=(a-b)^2+2b(a-b)$, we have
\begin{eqnarray}\label{ineq3}
\left\vert\phi_n(\mathcal{U}_{n,p})-\phi(\mathcal{U}_{n,p})\right\vert\
&\leq& \bigg\{\sum_{j=1}^{k}\left\vert\frac{n(1+(k-2)P_j)}{n_j-1}-\rho^{-1}_j-k+2\right\vert\nonumber\\
&&+\left(\sum_{j=1}^{k}\left\vert\frac{1-P_j}{\sqrt{P_j}}-\frac{1-\rho_j}{\sqrt{\rho_j}}\right\vert\right)^2\nonumber\\
& &+ 2\sum_{j=1}^{k} \sum_{\ell=1}^{k}\frac{1-\rho_\ell}{\sqrt{\rho_\ell}}\left\vert\frac{1-P_j}{\sqrt{P_j}}-\frac{1-\rho_j}{\sqrt{\rho_j}}\right\vert\nonumber\\
& &+\left(\sum_{j=1}^{k}\left\vert P_j^{-1/2}- \rho_j^{-1/2} \right\vert\right)^2\nonumber\\
& &+ 2\sum_{j=1}^{k} \sum_{\ell=1}^{k} \rho_\ell^{-1/2}\left\vert P_j^{-1/2}- \rho_j^{-1/2}\right\vert\nonumber\\
& &+\left(\sum_{j=1}^{k}\left\vert P_j^{1/2}- \rho_j^{1/2} \right\vert\right)^2\nonumber\\
& &+ 2\sum_{j=1}^{k} \sum_{\ell=1}^{k} \rho_\ell^{1/2}\left\vert P_j^{1/2}- \rho_j^{1/2}\right\vert\,\bigg\}\left\vert \mathcal{U}_{n,p}\right\vert^2
\end{eqnarray}
Since  $\mathcal{U}_{n,p}$ (resp. $\mathcal{V}_{n,p}$) converges in distribution (resp. in probability)  to   $\mathcal{U}_{p}$ (resp. ${\large \textrm{\textbf{1}}}_k$), we deduce from (\ref{ineq1}), (\ref{ineq2}), (\ref{ineq3}) and the continuity of $\phi$ and  $\psi$ that  $\widehat{S}_n^{(q)}-S_q$ converges in probability to $0$  as $\min_{1\leq j\leq k}(n_j)\rightarrow +\infty$ and, consequetly, that $\widehat{S}_n^{(q)}\stackrel{\mathscr{D}}{\rightarrow}S_q$  as $\min_{1\leq j\leq k}(n_j)\rightarrow +\infty$.  Therefore, there exists $N_1$ such that, for 
 $\min_{1\leq j\leq k}(n_j)\geq N_1$, one has
\begin{equation}\label{eps2}
\left\vert\varphi_{\widehat{S}_n^{(q)}}(t)-\varphi_{S_{(q)}}(t)\right\vert<\frac{\varepsilon}{3}.
\end{equation}
Thirdly,  let us consider  $S=\sum_{p=1}^{+\infty}\lambda_p\mathcal{W}_{p}$ and show that    $\widehat{S}^{(q)}\stackrel{\mathscr{D}}{\rightarrow}S_q$  as $q\rightarrow +\infty$. We have
\begin{eqnarray*}
\left\vert\varphi_{S_q}(t)-\varphi_{S}(t)\right\vert&\leq& 
\mathbb{E}\left(\left\vert e^{itS_q}-e^{itS}\right\vert\right)\leq\left\vert t\right\vert\mathbb{E}\left(\left\vert  S_q-S \right\vert\right)
\leq\left\vert t\right\vert\sum\limits_{p=q+1}^{+\infty}\lambda_p\mathbb{E}\left(\left\vert  \mathcal{W}_{p} \right\vert\right)
\end{eqnarray*}
and
\begin{eqnarray*}
\mathbb{E}\left(\left\vert  \mathcal{W}_{p} \right\vert\right)&=&\mathbb{E}\left(\left\vert  \phi(\mathcal{U}_{p})-\psi({\large \textrm{\textbf{1}}}_k) \right\vert\right)\\
&\leq&\sum_{j=1}^{k}\left(\rho_j^{-1}+k-2\right)\,\mathbb{E}\left(Y_{j,p}^2\right)+\mathbb{E}\left(\left[\sum_{j=1}^{k}\frac{1-\rho_j}{\sqrt{\rho_j}}Y_{j,p}\right]^2\right)\\
& &+\mathbb{E}\left(\left[\sum_{j=1}^{k}\rho_j^{-1/2}Y_{j,p}\right]^2\right)+\mathbb{E}\left(\left[\sum_{j=1}^{k}\rho_j^{1/2}Y_{j,p}\right]^2\right)\\
& &+\sum_{j=1}^{k}\left(\rho_j^{-1}+k-2\right).
\end{eqnarray*}
Since $\mathbb{E}\left(Y_{j,p}Y_{\ell,p}\right)=\delta_{j\ell}$, it follows

\begin{eqnarray*}
\mathbb{E}\left(\left[\sum_{j=1}^{k}\frac{1-\rho_j}{\sqrt{\rho_j}}Y_{j,p}\right]^2\right)&=&\sum_{j=1}^{k}\frac{(1-\rho_j)^2}{\rho_j}\mathbb{E}\left(Y_{j,p}^2\right)\\
& &+\sum_{j=1}^{k}\sum_{\underset{\ell\neq j}{\ell=1}}^{k}\frac{(1-\rho_j)(1-\rho_\ell)}{\sqrt{\rho_j\rho_\ell}}\mathbb{E}\left(Y_{j,p}Y_{\ell,p}\right)\\
&=&\sum_{j=1}^{k}\frac{(1-\rho_j)^2}{\rho_j},
\end{eqnarray*} 
\begin{eqnarray*}
\mathbb{E}\left(\left[\sum_{j=1}^{k}\rho_j^{-1/2}Y_{j,p}\right]^2\right)&=&\sum_{j=1}^{k} \rho_j^{-1} \mathbb{E}\left(Y_{j,p}^2\right) +\sum_{j=1}^{k}\sum_{\underset{\ell\neq j}{\ell=1}}^{k} \rho_j^{-1/2} \rho_\ell^{-1/2}\mathbb{E}\left(Y_{j,p}Y_{\ell,p}\right)\\
&=&\sum_{j=1}^{k} \rho_j^{-1} ,
\end{eqnarray*}
\begin{eqnarray*}
\mathbb{E}\left(\left[\sum_{j=1}^{k}\rho_j^{1/2}Y_{j,p}\right]^2\right)&=&\sum_{j=1}^{k} \rho_j \mathbb{E}\left(Y_{j,p}^2\right) +\sum_{j=1}^{k}\sum_{\underset{\ell\neq j}{\ell=1}}^{k} \rho_j^{1/2} \rho_\ell^{1/2}\mathbb{E}\left(Y_{j,p}Y_{\ell,p}\right)\\
&=&\sum_{j=1}^{k} \rho_j=1 .
\end{eqnarray*}
Thus
\begin{eqnarray*}
\mathbb{E}\left(\left\vert  \mathcal{W}_{p} \right\vert\right)
&\leq&\sum_{j=1}^{k}\bigg\{2\left(\rho_j^{-1}+k-2\right)+\frac{(1-\rho_j)^2+1+\rho_j}{\rho_j}\bigg\}
\end{eqnarray*}
and, therefore, 
\begin{eqnarray*}
\left\vert\varphi_{S_q}(t)-\varphi_{S}(t)\right\vert&\leq& 
\sum_{j=1}^{k}\bigg\{2\left(\rho_j^{-1}+k-2\right)+\frac{(1-\rho_j)^2+1+\rho_j}{\rho_j}\bigg\}\sum\limits_{p=q+1}^{+\infty}\lambda_p.
\end{eqnarray*}
Since $\sum_{p=1}^{+\infty}\lambda_p<+\infty$, we deduce  that $\lim_{q\rightarrow +\infty}\left\vert\varphi_{S_q}(t)-\varphi_{S}(t)\right\vert=0$. Then,   there exists $q_1\in\mathbb{N}^\ast$ such that 
\begin{equation}\label{eps3}
\left\vert\varphi_{\widehat{S}_n}(t)-\varphi_{\widehat{S}_n^{(q)}}(t)\right\vert<\frac{\varepsilon}{3}
\end{equation}
for $q\geq q_1$. 
Putting $u=\max(q_0,q_1)$, $N_0=\max(n_1^0,\cdots,n_k^0)$  and using (\ref{eps1}), (\ref{eps2}) and (\ref{eps3}), we deduce that, if    $\min_{1\leq j\leq k}(n_j)\geq \max(N_0,N_1)$ then
\begin{eqnarray*}
\left\vert\varphi_{\widehat{S}_n}(t)-\varphi_{S}(t)\right\vert\leq \left\vert\varphi_{\widehat{S}_n}(t)-\varphi_{\widehat{S}_n^{(u)}}(t)\right\vert
+\left\vert\varphi_{\widehat{S}_n^{(u)}}(t)-\varphi_{S_{u}}(t)\right\vert
+\left\vert\varphi_{S_{u}}(t)-\varphi_{S}(t)\right\vert <\varepsilon.
\end{eqnarray*}
This show that  $\widehat{S}_n\stackrel{\mathscr{D}}{\rightarrow}S$  as $\min_{1\leq j\leq k}(n_j)\rightarrow +\infty$, where
\begin{eqnarray*}
S&=&\sum\limits_{p=1}^{+\infty}\lambda_{p}\bigg\{\sum_{j=1}^{k}\left(\rho_j^{-1}+(k-2)\right)Y_{j,p}^2+2\sum_{j=1}^{k}Y_{j,p}^2+\left[\sum_{j=1}^{k}\frac{1-\rho_j}{\sqrt{\rho_j}}Y_{j,p}\right]^2\nonumber\\
& &-\left[\sum_{j=1}^{k}\rho_j^{-1/2}Y_{j,p}\right]^2-\left[\sum_{j=1}^{k}\rho_j^{1/2}Y_{j,p}\right]^2-\sum_{j=1}^{k}\left(\rho_j^{-1}+(k-2)\right)\bigg\}\\
&=&\sum_{p=1}^{+\infty}\lambda_p\left\{(k-2)(Z_p-k)+\sum_{j=1}^{k}\left\{\rho_j^{-1}(Y^2_{p,j}-1)-2\sum_{\underset{\ell\neq j}{\ell=1}}^k\rho_\ell^{1/2}\rho_j^{-1/2}Y_{p,j}Y_{p,\ell}\right\}\right\}
\end{eqnarray*}
and $Z_{p}=\sum_{j=1}^kY_{j,p}^2\leadsto\chi^2_k$. 
\hfill$\Box$
\begin{rmk}
With this theorem we recover Theorem 12 of Gretton et al (2012). Indeed, for  $k=2$ the random variable to which $n\widehat{T}_n$ converges in distribution is
\begin{eqnarray*}
S&=&\sum_{p=1}^{+\infty}\lambda_p \bigg(\rho_1^{-1}(Y^2_{1,p}-1)+\rho_2^{-1}(Y^2_{2,p}-1)-2(\rho_2^{1/2}\rho_1^{-1/2}+\rho_1^{1/2}\rho_2^{-1/2})Y_{p,1}Y_{p,2}\bigg)\\
&=&\sum_{p=1}^{+\infty}\lambda_p \bigg(\rho_1^{-1}Y^2_{1,p}+\rho_2^{-1}Y^2_{2,p}-\rho_1^{-1}-\rho_2^{-1}-2\rho_1^{1/2}\rho_2^{1/2}(\rho_1^{-1}+\rho_2^{-1})Y_{p,1}Y_{p,2}\bigg)
\end{eqnarray*}
Since 
\[
\rho_1^{-1}+\rho_2^{-1}=\frac{\rho_1+\rho_2}{\rho_1\rho_2}=\frac{1}{\rho_1\rho_2},
\]
we obtain
\begin{eqnarray*}
S&=&\sum_{p=1}^{+\infty}\lambda_p \bigg(\rho_1^{-1}Y^2_{1,p}+\rho_2^{-1}Y^2_{2,p}-(\rho_1\rho_2)^{-1}-2\rho_1^{-1/2}\rho_2^{-1/2}Y_{p,1}Y_{p,2}\bigg)\\
&=&\sum_{p=1}^{+\infty}\lambda_p \bigg\{\left(\rho_1^{-1/2}Y_{1,p}-\rho_2^{-1/2}Y_{2,p}\right)^2 -(\rho_1\rho_2)^{-1} \bigg\},
\end{eqnarray*}
what is the result in the aforementioned  Theorem 12.
\end{rmk}
\begin{rmk}
it is difficult, if not impossible, to use the result in Theorem 1    for the practical  implementation of the proposed test. So, one can use subsampling methods  (see, e.g., Politis et al (1999), Berg et al (2010)) for computing $p$-values in order to perfom this test by using the introduced test statistic.
\end{rmk}
\section{Monte carlo simulations}
\noindent In this section, the finite sample performance of the proposed test is evaluated through   Monte Carlo simulations and compared to tests introduced by Zhang and Wu (2007). These authors proposed three tests for the $k$-sample problem, based on statistics denoted by $Z_A$,  $Z_B$ and  $Z_C$ obtained from the likelihood-ratio test statistic,  and showed that these tests are more powerful than the classical Kolmogorov-Smirnov, Cram\'er-von Mises and Anderson-Darling $k$-sample tests. We  estimate the powers of  our test and the three  aforementioned  tests   in the case where  $(k = 3)$ , and considering the four following cases: 
\begin{eqnarray*}
&&\textrm{\textbf{Case 1:} } \mathbb{P}_1=N(3,1),\hspace*{0.2cm}\mathbb{P}_2=Gamma(3,1)\hspace*{0.2cm}\textrm{ and }\hspace*{0.2cm} \mathbb{P}_3=Gamma(6,2);\\
&&\textrm{\textbf{Case 2:} } \mathbb{P}_1=N(0,1),\hspace*{0.2cm}\mathbb{P}_2=N(0,2)\hspace*{0.2cm}\textrm{ and }\hspace*{0.2cm} \mathbb{P}_3=N(0,4);\\
&&\textrm{\textbf{Case 3:} } \mathbb{P}_1=Uniform(0,1),\hspace*{0.2cm}\mathbb{P}_2=Beta(1,1.5)\hspace*{0.2cm}\textrm{ and } \hspace*{0.2cm} \mathbb{P}_3=Beta(1.5,1);\\
&&\textrm{\textbf{Case 4:} } \mathbb{P}_1=N(0,1),\hspace*{0.2cm}\mathbb{P}_2=N(0.3,1)\hspace*{0.2cm}\textrm{ and }\hspace*{0.2cm} \mathbb{P}_3=N(0.6,1).
\end{eqnarray*}
For all tests we take  the significance level $\alpha= 0.05$ and the empirical power is computed over $100$ independent replications. For our test, we used the gaussian  kernel $K(x,y)=\exp[-2(x-y)^2]$, and since the asymtotic distribution given in Therorem \ref{loi} is hard to simulate, we computed  the sampling distributions of $n\widehat{T}_n$ under $\mathscr{H}_0$ in order to compute the corresponding $p$-values. The results are given in Figures 1 to 4  that   plot  the  empirical power versus the total sample size $n=n_1+n_2+n_3$. They show that our test outperforms the three tests of Zhang and Wu (2007) in all cases.
\begin{figure}[h!]
\centering
\includegraphics[width=8cm,height=8cm]{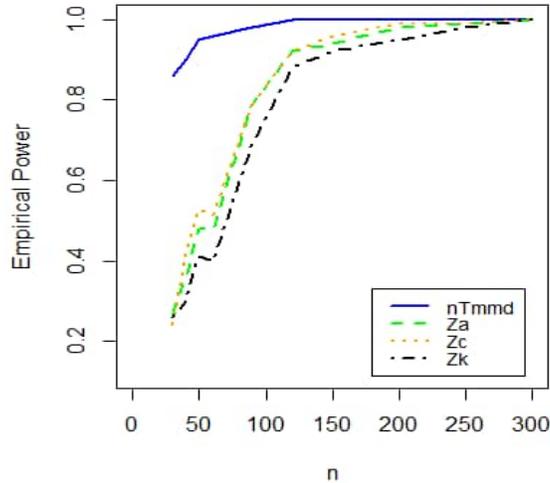}
\caption{Empirical power    versus $n$ for Case 1 with significance level $\alpha=0.05$}
\end{figure}
\begin{figure}[h!]
\centering
\includegraphics[width=8cm,height=8cm]{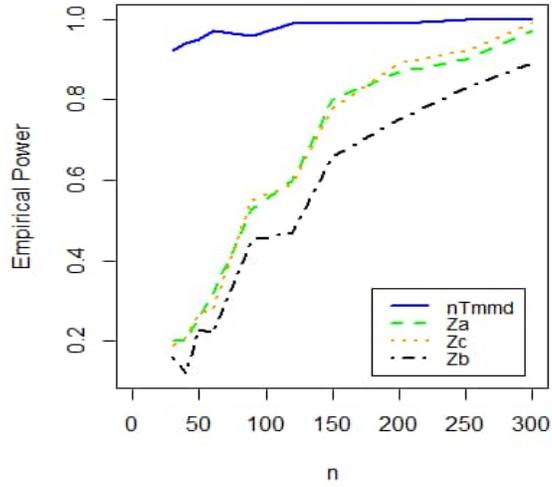}
\caption{Empirical power    versus $n$ for Case 2  with significance level $\alpha=0.05$}
\end{figure}
\begin{figure}[h!]
\centering
\includegraphics[width=8cm,height=8cm]{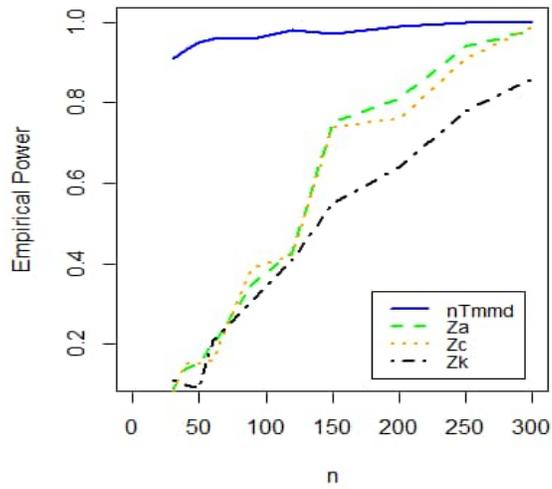}
\caption{Empirical power   versus $n$ for Case 3  with significance level $\alpha=0.05$}
\end{figure}
\begin{figure}[h!]
\centering
\includegraphics[width=8cm,height=8cm]{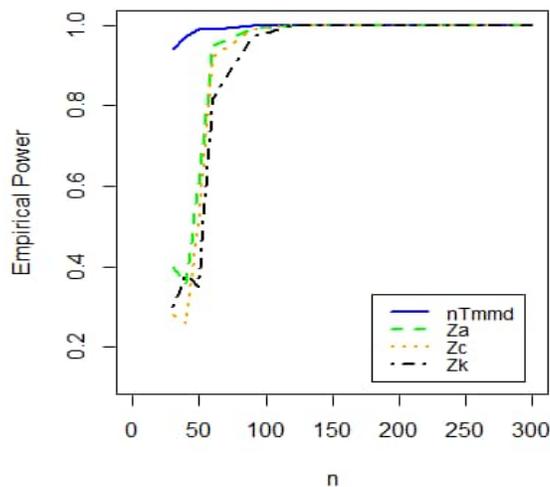}
\caption{Empirical power  versus $n$ for Case 4 with significance level $\alpha=0.05$}
\end{figure}


\section*{References}

\begin{thebibliography}{}


\bibitem{bA}Berlinet  A , Thomas-Agnan  C (2004)  Reproducing Kernel Hilbert Spaces in Probability and Statistics. Kluwer.

\bibitem{berg}Berg  A,  McMurry  TL, Politis  DN (2010)  Subsampling $p$-values. Stat. Probab. Lett. 80:1358-1364.

\bibitem{con}Conover  WJ (1980)  Practical nonparametric statistics. Wiley, New York.

\bibitem{gib}Gibbons  JD, Chakraborti S (2003)  Nonparametric statistical inference. Dekker, New York.

	\bibitem{GFT 12} Gretton A, Borgwardt  KM,  Rasch  MJ,  Sch$\ddot{\textrm{o}}$lkopf  B,  Smola  AJ (2012)  A kernel two-sample test. J.  Mach. Learn.  Res.  13:723-776.

     \bibitem{har}Harchaoui  Z,   Bach F,   Capp\'e  O , Moulines  E (2013)  Kernel-Based Methods for Hypothesis Testing: A Unified View. IEEE Signal Processing Magazine  30:87-97. 

\bibitem{kief}Kiefer  J (1959)  $k$-Sample analogues of the Kolmogorov-Smirnov, Cram\' er-von Mises tests. Ann. Math. Statist. 30:420-447.

	\bibitem{krus}Kruskal  WH, Wallis WA (1952)  Use of ranks in one-criterion analysis of variance. J. Amer. Statist. Assoc. 47:583-621.

\bibitem{pol}Politis  DN, Romano  JP, Wolf  M (1999)  Subsampling. Springer Verlag.

\bibitem{sch}Scholz FW, Stephens  MA (1987)  $k$-Sample Anderson-Darling tests. J. Amer. Statist. Assoc. 82:918-924.
	
\bibitem{wolf}Wolf  EH, Naus  JI (1973)  Tables of critical values for a $k$-sample   Kolmogorov-Smirnov  test  statistic. J. Amer. Statist. Assoc. 68:994-997.

          \bibitem{JZ 07}Zhang  J,  Wu Y (2007)  $ k$-Sample tests based on the likelihood ratio. Comput. Stat.   Data  Anal.  51:4682-4691.
\end{thebibliography}
\end{document}